# On (Co)homological Properties of Stone-Čech Compactifications of Completely Regular Spaces

Vladimer Baladze and Fridon Dumbadze

Department of Mathematics
of Batumi Shota Rustaveli State University, Batumi, Georgia

**Abstract**

Using the set of functionally open finite covers of completely regular spaces in the paper are constructed Čech type functional homology functor $\check{H}^F_q(-,-;G)\colon \boldsymbol{Top^2_{cr}} \to \boldsymbol{Ab}$ and functional cohomology functor $\hat{H}^q_F(-,-;G)\colon \boldsymbol{Top^2_{cr}} \to \boldsymbol{Ab}$ from the category of pairs of completely regular spaces and their completely closed subspaces to the category of abelian groups, defined Bokstein-Nowak type functional coefficient of cyclicity $\eta^F_G\colon \boldsymbol{Top_{cr}} \to \boldsymbol{N} \cup \{-1,\infty\}$ from the class of completely regular spaces to the set of integers $t \geq -1$, proved the equalities $\check{H}^F_n(X,A;G) = \check{H}_n(\beta X,\beta A;G)$, $\hat{H}^n_F(X,A;G) = \hat{H}^n(\beta X,\beta A;G)$and $\eta^F_G(X) = \eta_G(\beta X)$, where $\boldsymbol{N}$, $\check{H}_n(\beta X,\beta A;G)$, $\widehat{H}^n(\beta X,\beta A;G)$ and $\eta_G(\beta X)$ are the set of natural numbers, Čech homology group, Čech cohomology group and Bokstein-Nowak coefficient of cyclisity of Stone-Čech compactifications of pair $(X,A) \in ob(\boldsymbol{Top^2_{cr}})$ and space $X \in ob(\boldsymbol{Top_{cr}})$, respectively.

## 0. Introduction

The problems of compactification theory of topological spaces also lead to the necessity of the creation and development of such new (co)homological theories, whose methods effectively can be used in the study of specific problems of algebraic topology.

The investigation presented in this paper is centered around the following general
**Problem: Find necessary and sufficient conditions under which the spaces of the given class has the compactifications with the given topological properties.**

---

The first author was supported by grants FR/233/5-103/14 from Shota Rustaveli National Science Foundation (SRNSF) and Batumi Shota Rustaveli state University.

This problem for different topological properties was studied by several authors (cf. [A-N], [B$_1$], [B$_2$], [B-T], [Ba], [Bo], [E-S], [En$_1$], [En$_2$], [I], [K], [M-S], [Mi], [N], [P], [Sk], [Sm$_1$], [Sm$_2$], [Z], [V]).

The present paper is devoted to the study of this problem for property:

*n*-dimensional Čech homology (cohomology) group [E-S] and coefficient of cyclicity ([Bo],[N] ) of Stone- Čech compactification of completely regular space is a given group and integer , respectively.

A special aspect of this topic was considered in [E-S] , where was proved that if $(\tilde{X}, \tilde{A})$ is the Tychonoff compactification of closed pair $(X, A)$ of normal spaces, then

$$\check{H}_n^f(X, A; G) = \check{H}_n(\tilde{X}, \tilde{A}; G)$$

and

$$\hat{H}_f^n(X, A; G) = \hat{H}^n(\tilde{X}, \tilde{A}; G),$$

where $\check{H}_n^f(X, A; G)$ and $\hat{H}_f^n(X, A; G)$, $\hat{H}^n(\tilde{X}, \tilde{A}; G)$ and $\check{H}_n(\tilde{X}, \tilde{A}; G)$ respectively are n-dimensional Čech homology groups and n-dimensional Čech cohomology groups of pairs $(X, A)$ and $(\tilde{X}, \tilde{A})$, based on finite open covers. Also note that, in papers [B$_1$], [B-T] and [Mi] the obtained results include the characterizations of other homology and cohomology groups of extensions of spaces from some classes of spaces.

In this paper we define the Čech type functional homology and cohomology $\partial$-functor and $\delta$-cofunctor (see [E-S]) with coefficients in abelian group $\check{H}_n^F(-,-;G): \boldsymbol{Top}_{cr}^2 \to \boldsymbol{Ab}$ and $\hat{H}_F^n(-,-;G): \boldsymbol{Top}_{cr}^2 \to \boldsymbol{Ab}$ from the category $\boldsymbol{Top}_{cr}^2$ of pairs of completely regular spaces and their completely closed subspaces [Sm$_1$] to the category $\boldsymbol{Ab}$ of abelian groups.The definitions of these functors are based on the set of functionally open finite covers of pairs $(X, A) \in ob(\boldsymbol{Top}_{cr}^2)$.

The main results of the paper is the following

**Theorem 2.1** *For each pair $(X, A) \in ob(\boldsymbol{Top}_{cr}^2)$, one has*

$$\check{H}_n^F(X, A; G) = \check{H}_n(\beta X, \beta A; G)$$

*and*

$$\hat{H}_F^n(X, A; G) = \hat{H}^n(\beta X, \beta A; G).$$

In the paper also is defined the functional coefficient of cyclicity $\eta_G^F(X)$ of completely regular space $X$ and proved the following

**Theorem 2.2** *For each $X \in ob(\boldsymbol{Top}_{cr})$ holds the equality*

$$\eta_G^F(X) = \eta_G(\beta X).$$

Note that, here $\eta_G(\beta X)$ is coefficient of cyclisity of $\beta X$, defined and studied by M.F. Bokstein [Bo] and S. Novak [No].

Thus, the functional (co)homology groups intrinsically, in terms of functionaly open sets of completely regular spaces describe the Čech (co)homology groups and coefficients of cyclicity of Stone – Čech compactifications of completely regular spaces. In particular, a pair $(X, A) \in ob(\boldsymbol{Top}_{cr}^2)$ has the Stone – Čech compactification with Čech (co)homology group isomorphic to the abelian group $M$ if and only if $(X, A)$ has functional (co)homology group isomorphic to abelian group $M$. Besides, $X \in ob(\boldsymbol{Top}_{cr})$ has the Stone-Čech compactification with coefficient of cyclicity equal to the integer $n \geq -1$ if and only if $X$ has founctional coefficient of cyclicity equal to the integer $n \geq -1$.

In view of this results the relations $\check{H}_n^F(X,A;G) \approx M$, $\hat{H}_F^n(X,A;G) \approx M$ and $\eta_G^F(X) = n$ are the internal necessary and sufficient conditions on $(X,A)$ and $X$ so that (X,A) and X have Stone- Čech compactification with $\check{H}_n(\beta X, \beta A; G) \approx M$, $\hat{H}^n(\beta X, \beta A; G) \approx M$ and $\eta_G(\beta X) = n$.

Finally, let us note that without any further reference in the case of necessity we will use definitions and results from the books General Topology [En$_1$] and Algebraic Topology [E-S].

# 1. The Čech Type functional homology and cohomology functors

In this section of the paper using the methods developed in [E-S] and [En$_1$] we will give the construction of Čech Type functional homology and cohomology funtors on the category $(\boldsymbol{Top}_{cr}^2)$ and deduce some of their consequences.

Let $(X,A)$ be a pair consisting of completely regular space $X$ and its completely closed subset $A$. By $cov_F(X)$ denote the set of functionally open finite covers of space $X$. Let $\alpha = \{\alpha_v\}_{v\in V_\alpha} \in cov_F(X), |V_\alpha| < \chi_0$ and $s \subset V_\alpha$ be a subset of $V_\alpha$.. By $car_\alpha(s)$ denote the carrier of $s$. By definition, $car_\alpha(s) = \bigcap_{v\in s} \alpha_v$. Consider the nerve $X_\alpha$ of $\alpha$ and its subcomplex $A_\alpha$ consisting of all finite subsets $s \subset V_\alpha$ such that $car_\alpha(s) \cap A \neq \emptyset$ .

The chain and cochain groups $C_n(X_\alpha, A_\alpha; G)$ with coefficients in abelian group $G$ are defined in the usual way [E-S].

Now we assign to the simplicial pair $(X_\alpha, A_\alpha)$ homology and cohomology groups $H_n(X_\alpha, A_\alpha; G)$ and $H^n(X_\alpha, A_\alpha; G)$ with coefficients in group $G$ (see[E-S], ch. V1, 3.9).

Let $C_n(X_\alpha)$ be the free abelian group generated by all ordered $n$-simplices of $X_\alpha$ and $C_n(A_\alpha)$ be the free abelian subgroup generated by all ordered $n$-simplices of $A_\alpha$. By $C_n(X_\alpha, A_\alpha)$ denote the quotient group $C_n(X_\alpha)/C_n(A_\alpha)$. For each generator $(v_0, v_{1,\dots,}v_n)$ of $C_n(X_\alpha)$ is defined its boundary

$$\partial(v_0, v_{1,\dots,}v_n) = \sum_{i=0}^{n}(-1)^i\,(v_0, v_1, \dots, \check{v}_i, \dots v_n).$$

Hence, there exists boundary homomorphism $\partial_n: C_n(X_\alpha) \to C_{n-1}(X_\alpha)$ with properties $\partial_n(C_n(A_\alpha)) \subset C_{n-1}(A_\alpha)$ and $\partial_n(C_n(X_\alpha, A_\alpha)) \subset C_{n-1}(X_\alpha, A_\alpha)$. Thus, we have ordered chain complexes $C(X_\alpha) = \{C_n(X_\alpha), \partial_n\}$, $C(A_\alpha) = \{C_n(A_\alpha), \partial_n\}$ and $C(X_\alpha, A_\alpha) = \{C_n(X_\alpha, A_\alpha), \partial_n\}$.

Let

$$C_n(X_\alpha; G) = C_n(X_\alpha) \otimes G, C_n(A_\alpha; G) = C_n(A_\alpha) \otimes G$$

and

$$C_n(X_\alpha, A_\alpha; G) = C_n(X_\alpha, A_\alpha) \otimes G,$$

where symbol $\otimes$ denotes the tensor product of groups. Consequently, we have obtained the chain complexes

$$C(X_\alpha; G) = \{C_n(X_\alpha; G), \partial_n \otimes 1_G\},$$

$$C(A_\alpha;G)=\{C_n(A_\alpha;G),\partial_n\otimes 1_G\},$$

$$C(X_\alpha,A_\alpha;G)=\{C_n(X_\alpha,A_\alpha;G),\partial_n\otimes 1_G\}.$$

For simplicity the homomorphisms $\partial_n\otimes 1_G$ again denote by $\partial_n$.

As in ([E-S], ch.5,11) we can describe cochain complexes $\{C^n(X_\alpha;G),\delta^n\}$, $\{C^n(A_\alpha;G),\delta^n\}$ and $\{C^n(X_\alpha,A_\alpha;G),\delta^n\}$.

Let $H_n(X_\alpha,A_\alpha;G)=ker\partial_n/im\partial_{n+1}$ and $H^n(X_\alpha,A_\alpha;G)=ker\,\delta^n/Im\delta^{n-1}$. The boundary homomorphism $\partial: H_n(X_\alpha,A_\alpha;G)\to H_{n-1}(A_\alpha;G)$ and coboundary homomorphism $\delta: H^n(A_\alpha;G)\to H^{n+1}(X_\alpha,A_\alpha;G)$ are defined in the usual way. Also note that each simplicial map $f_\beta:(X_\alpha,A_\alpha)\to(Y_\beta,B_\beta)$ induces the homomorphism $f_{\beta *}: H_n(X_\alpha,A_\alpha;G)\to H_n(Y_\beta,B_\beta;G)$ and $f_\beta^*: H^n(Y_\beta,B_\beta;G)\to H^n(X_\alpha,A_\alpha;G)$.

Let $\alpha'\in cov_F(X)$ be a refinement of $\alpha\in cov_F(X)$ (notation $\alpha\le\alpha'$). Note that the set $cov_F(X)$ is a directed set with respect to the refinement relation $\alpha\le\alpha'$.

Any two refinement projection maps induce contiguous simplicial maps from pair $(X_{\alpha'},A_{\alpha'})$ to pair $(X_\alpha,A_\alpha)$. Consequently, any two refinement projection maps $p_{\alpha\alpha'}$ and $p'_{\alpha\alpha'}$ induce the unique homomorphisms

$$p_{\alpha\alpha'*}=p'_{\alpha\alpha'*}: H_n(X_{\alpha'},A_{\alpha'};G)\to H_n(X_\alpha,A_\alpha;G)$$

and

$$p^*_{\alpha\alpha'}=p'^*_{\alpha\alpha'}: H^n(X_\alpha,A_\alpha;G)\to H^n(X_{\alpha'},A_{\alpha'};G).$$

For each $\alpha\in cov_F(X)$ the homomorphisms $p_{\alpha\alpha*}$ and $p^*_{\alpha\alpha}$ are identity homomorphisms and for each triple $\alpha\le\alpha'\le\alpha''$, $p_{\alpha\alpha''*}=p_{\alpha\alpha'*}p_{\alpha'\alpha''*}$ and $p^*_{\alpha\alpha''}=p^*_{\alpha'\alpha''}\circ p^*_{\alpha\alpha'}$. Thus, there exists the inverse system $\left(H_n(X_\alpha,A_\alpha;G),p_{\alpha\alpha'*},cov_F(X)\right)$ and the direct system $(H^n(X_\alpha,A_\alpha;G),p^*_{\alpha\alpha'},cov_F(X))$ for pair $(X,A)\in ob(\boldsymbol{Top}^2_{cr})$.

We have the following

**Proposition 1.1.** *Let $A$ be a completely closed subset of space $X\in ob(\boldsymbol{Top}_{cr})$. The set of covers $\{\alpha_v\cap A|\alpha_v\in\alpha\in cov_F(X)\}$ is cofinal subset of the set $cov_F(A)$.*

**Proof.** Let $\alpha=\{\alpha_{v_i}\}_{v_i\in V_\alpha}, i=1,2,\dots,n$ be a functionally open cover of $A$. There exists a functionally closed shrinking $\beta=\{\beta_{v_i}\}_{v_i\in V_\alpha}$ of the cover $\alpha$ (see [En$_1$], Theorem 7.1.5). For each pair $\left(A\backslash\alpha_{v_i},\beta_{v_i}\right), v_i\in V_\alpha$ there exists a continuons map $f_{v_i}: A\to[0,1]$ such that $f_{v_i}\left(A\backslash\alpha_{v_i}\right)=\{0\}$ and $f_{v_i}\left(\beta_{v_i}\right)=\{1\}$. By condition of proposition each map $f_{v_i}$ has continuous extension $\tilde f_{v_i}: X\to[0,1]$. The family

$$\{\tilde f_{v_1}^{-1}\left(\left(\tfrac{1}{2},1\right]\right),\tilde f_{v_2}^{-1}\left(\left(\tfrac{1}{2},1\right]\right),\dots,\tilde f_{v_n}^{-1}\left(\left(\tfrac{1}{2},1\right]\right),\cap_{i=1}^n\tilde f_{v_i}^{-1}([0,1))\}$$

is a functionally open finite cover of $X$. Note that $A\cap\tilde f_{v_i}^{-1}\left(\left(\tfrac{1}{2},1\right]\right)\subset\alpha_{v_i}$ for each $i=1,2,\dots,n$ and $A\cap\cap_{i=1}^n\tilde f_{v_i}^{-1}\left([0,1)\right)=\emptyset$. Thus,

$$\{A\cap\tilde f_{v_i}^{-1}\left(\left(\tfrac{1}{2},1\right]\right),A\cap\tilde f_{v_2}^{-1}\left(\left(\tfrac{1}{2},1\right]\right),\dots,A\cap\tilde f_{v_n}^{-1}\left(\left(\tfrac{1}{2},1\right]\right),A\cap\cap_{i=1}^n\tilde f_{v_i}^{-1}([0,1))\}$$

is a functionally open finite cover of $A$ and it is a refinement of $\alpha$. □

The limit groups

$$\widetilde H_n^F(X,A;G)=\lim_{\leftarrow}\{H_n(X_\alpha,A_\alpha;G),p_{\alpha\alpha'*},cov_F(X)\}$$

and

$$\widehat{H}_F^n(X,A;G) = \lim_{\rightarrow}\{H^n(X_\alpha,A_\alpha;G), p^*_{\alpha\alpha'}, cov_F(X)\}$$

are called $n$-dimensional functional homology and cohomology groups of pair $(X,A) \in ob(\boldsymbol{Top^2_{cr}})$ with coefficients in abelian group $G$, respectively.

Let $(f:(X,A) \to (Y,B)) \in Mor_{\boldsymbol{Top^2_{cr}}}((X,A),(Y,B))$ and $\beta = \{\beta_v\}_{v\in V_\beta} \in cov_F(Y)$. The family $\alpha = \{\alpha_v\}_{v\in V_\alpha}$, where $\alpha_v = f^{-1}(\beta_v)$ for each $v \in V_\alpha = V_\beta$, is functionally open cover of $X$. Note that $X_\alpha$ is a subcomplex of $Y_\beta$ and $A_\alpha$ is a subcomplex of $B_\beta$. Hence, there exists the inclusion simplicial map $f_\beta:(X_\alpha,A_\alpha) \to (Y_\beta,B_\beta)$.

Let $\beta' \in cov_F(Y)$ and $\beta \le \beta'$. Then any refinement projection $p_{\beta\beta'}$ of $\beta'$into $\beta$ induces a refinement projection $p_{\alpha\alpha'}$ of $\alpha'$ into $\alpha$. From equality $f_\beta.p_{\alpha\alpha'} = p_{\beta\beta'}.f_{\beta'}$ follows that

$$f_{\beta*}.p_{\alpha\alpha'*} = p_{\beta\beta'*}.f_{\beta'*}$$

and

$$f_\beta{}^*.p^*{}_{\alpha\alpha'} = p_{\beta\beta'}{}^*.f^*{}_{\beta'}$$

Consequently, the map $f^{-1}: cov_f(Y) \to cov_F(X)$ and the homomorphisms

$$f_{\beta*}: H_n(X_\alpha,A_\alpha;G) \to H_n(Y_\beta,B_\beta;G)$$

and

$$f_\beta^*: H^n(Y_\beta,B_\beta;G) \to H^n(X_\alpha,A_\alpha;G)$$

form morphisms

$$(f_{\beta*},f^{-1}):(H_n(X_\alpha,A_\alpha;G),p_{\alpha\alpha'*},cov_F(X)) \to (H_n(Y_\beta,B_\beta;G),p_{\beta\beta'*},cov_F(Y))$$

and

$$(f_\beta^*,f^{-1}):(H^n(Y_\beta,B_\beta;G),p^*_{\beta\beta'},cov_F(Y)) \to (H^n(X_\alpha,A_\alpha;G),p^*_{\alpha\alpha'},cov_F(X)).$$

Let

$$f_* = \lim_{\leftarrow}(f_{\beta*},f^{-1}): \check{H}_n^F(X,A;G) \to \check{H}_n(Y,B;G)$$

and

$$f^* = \lim_{\rightarrow}(f_\beta^*,f^{-1}): \widehat{H}_F^n(Y,B;G) \to \widehat{H}_F^n(X,A;G).$$

The homomorphisms $f_*$ and $f^*$ are called the homomorphisms induced by the map $f:(X,A) \to (Y,B)$. Now we construct Čech type functional homology and cohomology founctors

$$\check{H}_n^F(-,-;G): \boldsymbol{Top^2_{cr}} \to \boldsymbol{Ab}$$

and

$$\widehat{H}_F^n(-,-;G): \boldsymbol{Top^2_{cr}} \to \boldsymbol{Ab}.$$

By definition,

$$\check{H}_n^F(-,-;G)((X,A)) = \check{H}_n^F(X,A;G), (X,A) \in ob(\boldsymbol{Top^2_{cr}}),$$
$$\check{H}_n^F(f) = f_*, f \in Mor_{\boldsymbol{Top^2_{cr}}}((X,A),(Y,B)),$$
$$\widehat{H}_F^n(-,-;G)((X,A)) = \widehat{H}_F^n(X,A;G), (X,A) \in ob(\boldsymbol{Top^2_{cr}}),$$
$$\widehat{H}_{cr}^n(f) = f^*, f \in Mor_{\boldsymbol{Top^2_{cr}}}((X,A),(Y,B)).$$

This implies that for each continuous maps $f \in Mor_{\boldsymbol{Top^2_{cr}}}((X,A),(Y,B))$, $g \in Mor_{\boldsymbol{Top^2_{cr}}}((Y,B),(Z,C))$ and $1_{(X,A)} \in Mor_{\boldsymbol{Top^2_{cr}}}((X,A),(X,A))$ hold the following equalities

$$\check{H}_n^F(g\circ f) = \check{H}_n^F(g)\circ\check{H}_n^F(f),$$
$$\widehat{H}_F^n(g\circ f) = \widehat{H}_F^n(f)\circ\widehat{H}_F^n(g),$$
$$\check{H}_n^F\big(1_{(X,A)}\big) = 1_{\check{H}_n^F(X,A;G)},$$
$$\widehat{H}_F^n\big(1_{(X,A)}\big) = 1_{\widehat{H}_F^n(X,A;G)}.$$

**Theorem 1.2**. *There exist functional homology and functional cohomology functors*

$$\check{H}_n^F(-,-;G)\colon \boldsymbol{Top_{cr}^2} \to \boldsymbol{Ab}$$

*and*

$$\widehat{H}_F^n(-,-;G)\colon \boldsymbol{Top_{cr}^2} \to \boldsymbol{Ab},$$

*respectively.*

**Proof**. The proof follows from the above given discussion and hence is omitted. □

Let $(X,A) \in ob(\boldsymbol{Top_{cr}^2})$. Consider two covers $\alpha,\alpha' \in cov_F(X)$ such that $\alpha \le \alpha'$. A refinement projection $p_{\alpha\alpha'}$ $of\ \alpha'$ into $\alpha$ induces the unique homomorphisms

$$p_{\alpha\alpha'*}\colon H_n(X_{\alpha'},A_{\alpha'};G) \to H_n(X_\alpha,A_\alpha;G)$$

and

$$p_{\alpha\alpha'}^*\colon H^n(X_\alpha,A_\alpha;G) \to H^n(X_{\alpha'},A_{\alpha'};G).$$

Let $V_\alpha^A = \{v \in V_\alpha | \alpha_v \cap A \neq \emptyset\}$. According to ([E-S], Ch. IX, 7) define a cover $\alpha' = \{\alpha_v \cap A\}_{v\in V_\alpha^A}$ of $A$ and an order preserving map $\varphi\colon cov_F(X) \to cov_F(A)$. By definition,$\varphi(\alpha) = \alpha'$. The correspondence $\alpha_v \to \varphi(\alpha)_v$ induces the identification of simplicial complexes $A_\alpha$ and $A_{\varphi(\alpha)}$. Let $\varphi_\alpha\colon A_{\varphi(\alpha)} \to A_\alpha$ be the inverse of this identification. If $\alpha \le \alpha'$ are covers of $X$, then the homomorphisms

$$p_{\alpha\alpha'*}\colon H_n\big(A_{\varphi(\alpha')};G\big) \to H_n\big(A_{\varphi(\alpha)};G\big),\ p_{\alpha\alpha'_*}\colon H_n(A_{\alpha'};G) \to H_n(A_\alpha;G)$$

and the homomorphisms

$$p_{\alpha\alpha'}^*\colon H^n\big(A_{\varphi(\alpha)};G\big) \to H^n\big(A_{\varphi(\alpha')};G\big),\ p_{\alpha\alpha'}^*\colon H^n(A_\alpha;G) \to H^n(A_{\alpha'};G)$$

are independent of the choice of refinement projection. Thus, hold the following equalities $p_{\alpha\alpha'*}\circ\varphi_{\alpha'*} = \varphi_{\alpha*}\circ p_{\varphi(\alpha)\varphi(\alpha')*}$ and $\varphi_{\alpha'}^*\circ p_{\alpha\alpha'}^* = p_{\varphi(\alpha)\varphi(\alpha')}^*\circ\varphi_\alpha^*$ for each pair $\alpha \le \alpha'$of $cov_F(X)$.

Let

$$\{H_n(A_\alpha;G), p_{\alpha\alpha'*}, cov_F(X)\}, \{H_n(A_\alpha), p_{\alpha\alpha'*}, cov_F(A)\}$$

and

$$\{H^n(A_\alpha;G),\ p_{\alpha\alpha}^*, cov_F(X)\}, \{H^n(A_\alpha;G),\ p_{\alpha\alpha}^*, cov_F(A)\}$$

be the inverse and direct systems, respectively. Their limits will be denoted by $\check{H}_n^F(A;G)_X, \check{H}_n^F(A;G)$ and $\widehat{H}_F^n(A;G)_X, \widehat{H}_F^n(A;G)$.

The map $\varphi$ and homomorphisms $\varphi_{\alpha*}$ and $\varphi_\alpha^*$ induce the maps of inverse and direct systems $(\varphi_{\alpha*},\varphi)$ and $(\varphi_\alpha^*,\varphi)$, respectively. It is clear that the limit homomorphisms

$$\Phi_n = \lim_{\leftarrow}(\varphi_{\alpha*},\varphi) : \check{H}_n^F(A,G) \to \check{H}_n^F(A;G)_X$$

and

$$\Phi^n = \lim_{\rightarrow}\ (\varphi_\alpha^*,\varphi) : \widehat{H}_F^n(A;G)_X \to \widehat{H}_F^n(A,G)$$

are isomorphisms.

Let $\alpha, \alpha' \in cov_F(X)$ and $\alpha \leq \alpha'$. Each refinement projection $p_{\alpha\alpha'}$ of $\alpha'$ to $\alpha$ induces the following commutative diagram of chain complexes

$$\begin{array}{ccccccccc}
0 \to & C(A_{\alpha'}; G) & \to & C(X_{\alpha'}; G) & \to & C(X_{\alpha'}, A_{\alpha'}; G) & \to 0 \\
& \downarrow p_{\alpha\alpha'} & & \downarrow p_{\alpha\alpha'} & & \downarrow p_{\alpha\alpha'} & \\
0 \to & C(A_\alpha; G) & \to & C(X_\alpha; G) & \to & C(X_\alpha, A_\alpha; G) & \to 0,
\end{array}$$

where short exact rows are induced by simplicial inclusions $A_\alpha \to X_\alpha$, $X_\alpha \to (X_\alpha, A_\alpha)$, $A_{\alpha'} \to X_{\alpha'}$ and $X_{\alpha'} \to (X_{\alpha'}, A_{\alpha'})$.

The commutativity relation holds in the following rectangle

$$\begin{array}{ccc}
H_n(X_{\alpha'}, A_{\alpha'}; G) & \xrightarrow{\partial_{\alpha'}} & H_{n-1}(A_{\alpha'}; G) \\
\downarrow p_{\alpha\alpha'*} & & \downarrow p_{\alpha\alpha'*} \\
H_n(X_\alpha, A_\alpha; G) & \xrightarrow{\partial_\alpha} & H_{n-1}(A_\alpha; G).
\end{array}$$

Hence, the family $\{\partial_\alpha\}$ induces a morphism of inverse systems and a limit homomorphism $\partial_{n,X}: \check{H}_n^F(X, A; G) \to \check{H}_{n-1}^F(A; G)_X$. Let

$$\partial = \Phi_n^{-1} \circ \partial_{n,X}: \check{H}_n^F(X, A; G) \to \check{H}_{n-1}^F(A; G).$$

Analogously, we can define a homomorphism $\delta_X^{n-1}: \widehat{H}_F^{n-1}(A; G)_X \to \widehat{H}_F^n(X, A; G)$. Let assume that

$$\delta = \delta_X^{n-1} \circ (\Phi^{n-1})^{-1}: \widehat{H}_F^{n-1}(A; G) \to \widehat{H}_F^n(X, A; G).$$

Then we have the following

**Theorem 1.3.** For every map $f: (X, A) \to (Y, B)$, integer $n$ and abelian group $G$ the following rectangles

$$\begin{array}{ccc}
\check{H}_n^F(X, A; G) & \xrightarrow{\partial} & \check{H}_{n-1}^F(A; G) \\
f_* \downarrow & & \downarrow (f_{|A})_* \\
\check{H}_n^F(Y, B; G) & \xrightarrow{\partial} & \check{H}_{n-1}^F(B; G)
\end{array}$$

and

$$\begin{array}{ccc}
\widehat{H}_F^{n-1}(B; G) & \xrightarrow{\delta} & \widehat{H}_F^{n-1}(Y, B; G) \\
(f_{|A})^* \downarrow & & \downarrow f^* \\
\widehat{H}_F^{n-1}(A; G) & \xrightarrow{\delta} & \widehat{H}_F^{n-1}(X, A; G)
\end{array}$$

are commutative.

**Proof:** The proof follows from above given discussion and hence is omitted. □

Therefore the functional homology and cohomology functors

$$\check{H}_n^F(-,-;G): \boldsymbol{Top_{cr}^2} \to \boldsymbol{Ab}$$

*and*

$$\widehat{H}_F^n(-,-;G): \boldsymbol{Top_{cr}^2} \to \boldsymbol{Ab}$$

are $\partial$ -functor and $\delta$-cofunctor in the sence of [E-S], respectively. This is a direct consequence of Theorems 1.2,1.3 and 1.4.

For each functionally open cover $\alpha$ of $X$ the inclusion maps $i: A \to X$ and $j: X \to (X, A)$ induce a short exact sequence of chain complexes

$$0 \to C(A_\alpha; G) \xrightarrow{i} C(X_\alpha; G) \xrightarrow{j} C(X_\alpha, A_\alpha; G) \to 0.$$

Consequently, there exist the long exact sequences

$$\dots \to H_n(A_\alpha; G) \xrightarrow{i_*} H_n(X_\alpha; G) \xrightarrow{j_*} H_n(X_\alpha, A_\alpha; G) \xrightarrow{\partial_*} H_{n-1}(A_\alpha; G) \to \cdots$$

and

$$\dots \to H^{n-1}(A_\alpha; G) \xrightarrow{\delta^*} H^n(X_\alpha, A_\alpha; G) \xrightarrow{j^*} H^n(X_\alpha; G) \xrightarrow{i^*} H^n(A_\alpha; G) \to \cdots \quad .$$

The limits of these inverse and direct systems

$$\dots \to \check{H}_n^F(A; G)_X \xrightarrow{i_*} \check{H}_n^F(X; G) \xrightarrow{j_*} \check{H}_n^F(X, A; G) \xrightarrow{\partial_*} \check{H}_{n-1}^F(A; G)_X \to \cdots$$

and

$$\dots \to \widehat{H}_F^{n-1}(A; G)_X \xrightarrow{\delta^*} \widehat{H}_F^n(X, A; G) \xrightarrow{j^*} \widehat{H}_F^n(X; G) \xrightarrow{i^*} \widehat{H}_F^n(A; G)_X \to \cdots$$

are of order 2 and exact , respectively. Therefore, we have the following

**Theorem 1.4.** *For any pair $(X, A) \in ob(\boldsymbol{Top_{cr}^2})$ the homology sequence*

$$\dots \to \check{H}_n^F(A; G) \to \check{H}_n^F(X; G) \to \check{H}_n^F(X, A; G) \to \check{H}_{n-1}^F(A; G) \to \cdots$$

*is of order 2 and the cohomology sequence*

$$\dots \to \widehat{H}_F^{n-1}(A; G) \to \widehat{H}_F^n(X, A; G) \to \widehat{H}_F^n(X; G) \to \widehat{H}_F^n(A; G) \to \cdots$$

*is exact.*

**Proof.** The proof follows from the above given discussion and hence is ommited. □

Let $(X, A, B)$ be a triple of spaces, i.e. $(X, A, B)$ consists of completely regular space $X$ and its completely closed subsets $A$ and $B$ with $A \supset B$.

The inclusion maps $i: A \to X$ and $j: X \to (X, A)$, $i': B \to X$ and $j': X \to (X, B)$, $i'': B \to A$ and $j'': A \to (A, B)$ induce the inclusion maps $\bar{\iota}: (A, B) \to (X, B)$ and $\bar{j}: (X, B) \to (X, A)$.For every integer $n$ consider the composition

$$\bar{\delta} = \delta^\circ j''^*: \widehat{H}_F^{n-1}(A, B; G) \to \widehat{H}_F^n(X, A; G).$$

There exists the Čech functional cohomology sequence of triple $(X, A, B)$. In particular, we have the following

**Theorem 1.5.** *The Čech functional cohomology sequence*

$$\dots \to \widehat{H}_F^{n-1}(A, B; G) \xrightarrow{\bar{\delta}} \widehat{H}_F^n(X, A; G) \xrightarrow{\bar{j}^*} \widehat{H}_F^n(X, B; G) \xrightarrow{\bar{\iota}^*} \widehat{H}_F^n(A, B; G) \to \cdots$$

*of triple $(X, A, B)$ is exact.*

**Proof.** The proof is dual and analogical to that of ([E-S], ch.I,10) and hence is ommited. □

Now, as an immediate consequence, we calculate the Čech functional (co)homology group of single points.

**Theorem 1.6**. *For the distinguished singleton space* $X = \{*\}$

$$\check{H}^F_n(X;G) = \begin{cases} 0, & n \neq 0 \\ G, & n = 0 \end{cases}$$

*and*

$$\hat{H}^n_F(X;G) = \begin{cases} 0, & n \neq 0 \\ G, & n = 0. \end{cases}$$

**Proof.** Let $\alpha$ be the cover of $X$ consisting of the single point . The inverse system of nerves of functionally open covers of $X$ has a single term $X_\alpha = N(\alpha)$ consisting of a single vertex. Consequently, $\check{H}^F_n(X;G) = H_n(N(\alpha);G)$ and $\hat{H}^n_F(X;G) = H^n(N(\alpha);G)$. Hence, $\check{H}^F_0(X;G) \approx G$ and $\check{H}^F_n(X;G) = 0$ for $n \neq 0$. Analogously, $\hat{H}^0_F(X;G) \approx G$ and $\hat{H}^n_F(X;G)$=0 for $n \neq 0$.□

## 2. The characterizations of Čech homology and cohomology groups and coefficients of cyclicity of Stone – Čech compactifications

Now we are mainly interested in the following question: How the Čech (co)homology groups and coefficients of cyclicity of Stone- Čech compactifications of completely regular spaces can be characterized intrinsically, in terms of functionally open subsets of completely regular spaces. The main result about the connection between Čech (co)homology groups of Stone – Čech compactifications of completely regular spaces and Čech functional (co)homology groups of completely regular spaces is the following

**Theorem 2.1** *Let* $(X, A) \in ob(\boldsymbol{Top^2_{cr}})$ *and* $(\beta X, \beta A)$ *be the pair of Stone- Čech compactifications of* $X$ *and* $A$*, respectively. Then*

$$\check{H}^F_n(X, A; G) = \check{H}_n(\beta X, \beta A; G)$$

*and*

$$\hat{H}^n_F(X, A; G) = \hat{H}^n(\beta X, \beta A; G).$$

**Proof.** Let $\alpha = \{\alpha_{v_i}\}_{v_i \in V_\alpha}$ be a functionally open finite cover of $X$. As in the proof of Proposition 1.1 we can show that there exists the functionally open finite cover of Stone- Čech compactification $\beta X$

$$\tilde{\alpha} = \{\tilde{f}^{-1}_{v_i}\left(\left(\tfrac{1}{2}, 1\right]\right), \cap_{i=1}^n \tilde{f}^{-1}_{v_i}([0,1))\}_{v_i \in V_\alpha}, \; i = 1, 2, \ldots, n$$

such that $\tilde{f}^{-1}_{v_i}\left(\left(\frac{1}{2}, 1\right]\right) \cap X \subset \alpha_{v_i}$ for each $i$= 1, 2,…,n and $(\cap_{i=1}^n \tilde{f}^{-1}_{v_i}([0,1)) \cap X = \emptyset$. Since $X$ is the dense subspace of $\beta X$, it follows that $\cap_{i=1}^n \tilde{f}^{-1}_{v_i}([0,1)) = \emptyset$. However, $\tilde{\alpha} \wedge X = \{\tilde{f}^{-1}_{v_i}\left(\left(\frac{1}{2}, 1\right]\right) \cap X\}_{v_i \in V_\alpha}$ is a functionally open finite cover of $X$ and it is a refinement of $\alpha$. Hence, for each $\alpha \in cov_F(X)$ there exists a cover $\tilde{\alpha} \in cov_F(\beta X)$ such that $\tilde{\alpha} \wedge X \geq \alpha$.

Let $\alpha = \{\alpha_{v_i}\}_{v_i \in V_\alpha}, i = 1,2,\dots,n$ be a finite open cover of $\beta X$. By Theorem 7.1.5 of [En$_1$] there exists a functionally open finite cover $\beta = \{\beta_{v_i}\}_{v_i \in V_\beta}$ of $\beta X$, which is a functionally open shrinking of $\alpha$ such that $\bar{\beta}_{v_i} \subset \alpha_{v_i}$ for each $i = 1,2,\dots,n$. Hence, $cov_F(\beta X)$ is cofinal subset of set $cov(\beta X)$. The set $\{\alpha \wedge X | \alpha \in cov_F(\beta X)\}$ is cofinal subset of set $cov_F(X)$. We have

$$\check{H}_n(\beta X, \beta A; G) = \lim_{\leftarrow}\{H_n((\beta X)_\beta, (\beta A)_\beta; G), p_{\beta\beta'_*}, cov_F(\beta X)\}$$

and

$$\check{H}_n^F(X, A; G) = \lim_{\leftarrow}\{H_n(X_{\beta\wedge X}, A_{\beta\wedge X}; G), p_{\beta\wedge X, \beta'\wedge X}, cov_F(\beta X)\}.$$

It is clear that the nerves $(\beta X)_\beta$ and $X_{\beta\wedge X}$ and $(\beta A)_\beta$ and $A_{\beta\wedge X}$ are isomorphic. These simplicial isomorphisms yield the isomorphisms

$$H_n((\beta X)_\beta, (\beta A)_\beta; G) \approx H_n(X_{\beta\wedge X}, A_{\beta\wedge X}; G), \beta \in cov_F(\beta X).$$

Applying the last isomorphisms we infer that $\check{H}_n^F(X, A; G) = \check{H}_n(\beta X, \beta A; G)$.
We analogously verify $\widehat{H}_F^n(X, A; G) = \widehat{H}^n(\beta X, \beta A; G)$.□

**Corollary 2.2.** *The Stone – Čech compactification* $(\beta X . \beta A)$ *of pair* $(X, A) \in ob(\boldsymbol{Top}_{cr}^2)$ *has Čech (co)homology group isomorphic to abelian group* $M$ *if and only if* $(X, A)$ *has Čech functional (co)homology group isomorphic to* $M$.

Now give the following

**Definition 2.3.** Let $G$ be an abelian group. A functional coefficient of cyclisity $\eta_G^F(X)$ of $X \in ob(\boldsymbol{Top}_{cr})$ with coefficients group $G$ is equal to $n$, if $\widehat{H}_F^m(X; G) = 0$ for each $m > n$ and $\widehat{H}_F^n(X; G) \neq 0$. For empty set $\emptyset$ we put $\eta_G^F(\emptyset) = -1$ and $\eta_G^F(X) = \infty$ if $X \neq \emptyset$ and for every $m$ there is an $n \geq m$ with $\widehat{H}_F^n(X; G) \neq 0$.

Thus, the functional coefficient of cyclicity of $X$ with coefficients group $G$ is a function $\eta_G^F: \boldsymbol{Top}_{cr} \to \boldsymbol{N} \cup \{-1, \infty\}: X \to n$, such that $\eta_G^F(X) = n$.

The following theorem is a direct consequence of Theorem 2.1.

**Theorem 2.4.** *For each* $X \in ob(\boldsymbol{Top}_{cr})$ *holds equality*

$$\eta_G^F(X) = \eta_G(\beta X).$$

**Proof.** Assume that $\eta_G(\beta X) = n$. Then for every $m > n$ and $\widehat{H}^m(\beta X; G) = 0$ and $\widehat{H}^n(\beta X; G) \neq 0$. From the relation $\widehat{H}^k(\beta X; G) = \widehat{H}_F^k(X; G)$ it follows that $\widehat{H}_F^m(X; G) = 0$ for every $m > n$ and $\widehat{H}_F^n(X; G) \neq 0$. Thus, $\eta_G^F(X) = n$. Consequently, $\eta_G^F(X) = \eta_G(\beta X)$.□

**Corollary 2.5.** *For each* $X \in ob(\boldsymbol{Top}_{cr})$ *the coefficient of cyclicity* $\eta_G(\beta X)$ *is equal to natural number* $n \geq -1$ *if and only if the functional coefficient of cyclicity* $\eta_G^F(X)$ *is equal to* $n$.

**Remark 2.6.** The construction of meaningfull cohomological dimension theory for completely regular spaces is one of main problems of modern topology. It is known that the theorems of classical cohomological dimension theory ([D], [H-W], [N], [Na]) of "good" (metric compact or Hausdorff compact) spaces often are impossible to generalize for completely regular spaces. In our next paper we will investigate this problem, using Čech functional cohomology theory. We define the small and large functional cohomological dimensions and establish their main properties.